\documentclass{elsart}

\usepackage{natbib}
\usepackage{amsfonts,amsmath}
\usepackage[dvips]{graphicx}

\newcommand{\nn}{\mathbf{n}}
\newcommand{\ttt}{\mathbf{t}}
\newcommand{\bx}{{\mathbf x}}
\newcommand{\bd}{{\mathbf d}}

\newcommand{\as}{{$\ast$}}

\begin{document}

\begin{frontmatter}
\title{Polyhedral Conditions for the Nonexistence of the MLE  for Hierarchical Log-linear Models}

\author[Berkeley]{Nicholas Eriksson}
\ead{eriksson@math.berkeley.edu}
\author[CMU,CMU2]{Stephen E. Fienberg}
\ead{fienberg@stat.cmu.edu}
\author[CMU]{Alessandro Rinaldo}
\ead{arinaldo@stat.cmu.edu}
\author[Berkeley]{Seth Sullivant}
\ead{seths@math.berkeley.edu}
\address[Berkeley]{Department of Mathematics, University of
  California, Berkeley }
\address[CMU]{Department of Statistics, Carnegie Mellon University}
\address[CMU2]{Center for Automated Learning and Discovery and Center for Computer and Communication Security, Carnegie Mellon University} 

\begin{abstract}
We provide a polyhedral description of the conditions for the existence of the
maximum likelihood estimate (MLE) for a hierarchical log-linear model. 
The MLE exists if and only if the
observed margins lie in the relative interior of the marginal cone.
Using this description, we give an algorithm for determining if the
MLE exists.  If the tree width is bounded, the
algorithm runs in polynomial time.  We also perform a computational
study of the case of three random variables under the no three-factor
effect model.
\end{abstract}


\end{frontmatter}

\section{Introduction}
 
In the analysis of contingency tables using log-linear models, 
the maximum likelihood estimate  (MLE) of  the underlying parameters 
(or equivalently of the expectations of the cell counts) 
plays a fundamental role for  computation, the assessment of
model fit, and model interpretation.  In particular, the existence of
the MLE is 
crucial for the determination of degrees of freedom  of  traditional
$\chi^2$ large sample 
approximations \citep[see, for example,][]{BFH75} and for exact or
approximate techniques for computing $p$-values.  If the MLE
does not exist, then the standard procedures and their approximations require alteration.

The characterizations of the conditions for the existence of the MLE
developed in the statistical literature are non-constructive, in the
sense that they do not directly lead to a numerical implementation
\cite[see][Appendix~B]{HAB74}.   As a result, the possibility of
the nonexistence of the MLE is rarely considered by practitioners and the
only available indication of it is a lack of convergence of the iterative
algorithms used to approximate the MLE.  

The problem of nonexistence has long been known to relate to the
presence of zero cell counts in the table, e.g., see  \cite{FIE70,HAB74,BFH75}.  
Zero counts arise frequently in large
sparse tables where the total sample size is small relative to the
number of cells in the table, e.g., see \cite{KOE86}.   Thus for small
contingency tables with a large sample size, the nonexistence 
of the MLE is a relatively infrequent problem.  This is because
for small contingency tables (nearly) all of the cell entries in the table
will be positive, which, as we will see, guarantees the existence of
the MLE.  However, the nonexistence of the MLE is a potentially
common problem in applications in the biological, medical, and social
sciences, where the 
contingency tables which arise are large and sparse.  Unfortunately,
in many such applications researchers ``collapse" large sparse tables
to form one of smaller dimension and/or size.  As \cite{BFH75} and
\cite{LAU96} make clear, such collapsing can lead to erroneous
statistical inferences about associations among the variables
displayed in the table. 

The goals of this paper are two-fold. First, we show that the 
nonexistence of the MLE is equivalent to the
margins of the observed contingency table lying on a facet of the
marginal cone of the  underlying hierarchical log-linear model.  
This polyhedral reinterpretation of
the problem immediately leads to easily implementable algorithms for
determining whether or not the MLE exists given an observed
contingency table and, in event the MLE does not exist, for
identifying those zero cell counts that cause the non-existence
problem.  We discuss these algorithms in Section 
\ref{sec:exist}.  From the practical standpoint, this characterization
gives a simple 
way to check whether or not the MLE exists before using  numerical
methods to estimate the MLE. 

The second goal of this paper is to alert the mathematical reader to a
rich source of combinatorial problems that arise from statistical
applications.  The polyhedral cones we are concerned with have
received attention in various guises (e.g., the ``correlation
polytope'' in \citet{deza-laurent} and the ``marginal polytope'' in
\citet{wainwright}).   Thus our particular problem of
deciding if a point in this cone is on a facet is a new variation on an
old theme.  Given recent computational advances, this also suggests
the problem of developing efficient algorithms for computing the
convex hulls of highly symmetric polyhedra.  We discuss these issues 
in Section \ref{sec:three}.

The outline for this paper is as follows.  In Section~\ref{sec:hier} we define
hierarchical models and the MLE, and we show that the MLE 
exists if and only if the observed margins belong to the relative
interior of a polyhedron. In Section~\ref{sec:exist} we use this fact
to describe an algorithm for checking the existence of the MLE.  The
algorithm uses linear programming and runs in polynomial time if the
tree width of the model is bounded.  Section~\ref{sec:three} focuses
on the study of the complexity of the problem for 3-way
tables. In particular, we consider the collapsing operation that preserves 
some combinatorial properties of a contingency table.

\section{Hierarchical models and the MLE}\label{sec:hier}

In this section, we introduce hierarchical models and the maximum
likelihood estimate and we show that the maximum likelihood estimate
exists if and only if certain polyhedral conditions are satisfied.
For this and the remaining sections, we assume the reader is familiar
with the basics of polyhedral geometry.  Two standard references are
\citet{ziegler} for basics on polyhedra and \citet{SCH98} for
algorithmic aspects including linear programming.  Our polyhedral
condition is a reformulation of a result of \citet{HAB74}.

Contingency tables are collections of non-negative integers arising
from cross-classifying a set of objects into  categories or cells
indexed by a set of labels $d$ corresponding to variables of interest 
\citep[see][]{BFH75,LAU96}.
More precisely, we get a $K$-way contingency table ${\bf n}$
by taking a sample of independent and identically distributed observations 
on  a vector of $K$ discrete random variables 
$(X_1, \ldots, X_K) $. The $j$th random variable $X_j$ takes values in the
set $[d_j] := \{1,2 \ldots, d_j\}$.  We call the various states of the random variables
\emph{levels}.   Let $d = \bigotimes_{j=1}^K [d_j]$. Thus each $i \in d$ identifies the
number ${\bf n}(i)$.  

Although the entries in the table ${\bf n}$ are integer-valued, we treat ${\bf n}$ 
 as an element of $\mathbb{R} ^ {d}$, the space of all real valued
functions on the multi-index set $d$ endowed with the usual inner
product
${\bf x}^{T}{\bf y} = \sum_{i \in d} {\bf x}(i) {\bf y}(i)$ for
${\bf  x},{\bf y} \in \mathbb{R} ^ {d}.$
For the remainder of the paper, we assume that the index set $d$
is linearized in some fashion, so that we can represent the table ${\bf n}$
 as a vector.

The statistical analysis of tables using log-linear models focuses on inference 
about parameters in a model or equivalently on 
 inferences about  the mean vector ${\bf m} =
\mathbb{E}({\bf n})$ of the observed table under the assumption that
${\bf m} > {\bf 0}$, so that $\mathbf{\mu} = \log \mathbf{m}$ is well
defined.  There are interesting extensions of the ideas in this paper to situations 
where we know a priori that some entries of ${\bf m}$are zero (c.f., \cite{BFH75,HAB74,FIE70}).

Log-linear models arise from assuming $\mu \in
\mathcal{M}$, where $\mathcal{M}$ is a $p$-dimensional linear subspace
$\mathcal{M} \subseteq \mathbb{R} ^ d$ such that ${\bf 1}_d \in
\mathcal{M}$. 
A common way of obtaining $\mathcal{M}$ is by specifying a
hierarchical model.  A hierarchical model is determined by a simplicial
complex $\Delta$ on $K$ vertices from which a 0-1 matrix $A_{\Delta}$
is constructed whose rows span $\mathcal{M}$ in the following way.   Let $\{
\mathcal{F}_1,\ldots,\mathcal{F}_f\}$ be the facets of $\Delta$ and,
for each $\mathcal{F}_s$ and $ i \in d$, let $d_{\mathcal{F}_s} =
\bigotimes_{j \in {\mathcal{F}_s}} [d_j]$ and $i_{\mathcal{F}_s}$ be the
restriction of $i$ to $d_{\mathcal{F}_s}$. Let $F_{\mathcal{F}_s}$ be
the set of functions on $d$ that depends on $i$ only
through $i_{\mathcal{F}_s}$.  That is,
\[
F_{\mathcal{F}_s} := 
\{ \bx \in \Rset^d \mid \bx(i) = \bx(j) \text{ for all } i,j \text{ with }  i_{\mathcal{F}_s}=j_{\mathcal{F}_s} \}.
\]
Then the linear subspace $\mathcal{M}$ corresponding to the
hierarchical log-linear model  $\Delta$ takes the form
\[
\mathcal{M}_{\Delta} = \sum_{\mathcal{F}_s \in \Delta} F_{\mathcal{F}_s}.
\]
Let $A_{\Delta}$ be a 0-1 matrix having dimension $v \times |d|$, where
$v = \sum_s \prod_{j \in \mathcal{F}_s} d_j$ and $|d|$ is the
cardinality of this index set.  Each row of $A_{\Delta}$ is indexed
by the pair $(\mathcal{F}_s, i_{\mathcal{F}_s})$ and is equal to the
indicator  function
$\chi(i_{\mathcal{F}_s})$, a  vector in $\mathbb{R}^d$ which
is 1 on coordinates $i_{\mathcal{F}_f}$ and 0 otherwise.  Then the
rows of $A_{\Delta}$ span  $\mathcal{M}_\Delta$,
so a hierarchical model can be identified by a collection of $K$
levels  ${\bf d} = (d_1,\ldots,d_K)$ and a simplicial complex $\Delta$
on $K$ nodes.

Data displayed in the form of contingency tables arise from various 
sampling  schemes involving the observations on the random variables 
\citep[see][]{BFH75,HAB74}. The results that  follow are  valid for 
the following three schemes:
\begin{description}
\item[Poisson Sampling.]
The total number $n= |{\bf n}|$ of counts is
random, where, for a non-negative vector ${\bf x}$, $|{\bf x}| =
\sum_i {\bf x}(i)$, and the counts are in fact independent Poisson random variables.
\item[Multinomial sampling.] The total number  $n= |{\bf n}|$ of counts is fixed by design.
\item[Product Multinomial sampling.]
Let $\mathcal{B} \subset \{1,\ldots,n \}$ and $d_{\mathcal{B}} =
\bigotimes_{j \in \mathcal{B}} d_j$, as above. For each $b \in
\mathcal{B}$, the number of counts $|{\bf n}(i_b)|$ is fixed by
design. Here, we assume, as is commonly done in the statistical literature,
 that $\mathcal{B}$ is always a face of
$\Delta$.
\end{description}

Given a table ${\bf n}$  on the fixed set of
levels ${\bf d} = (d_1,\ldots,d_K)$ and a simplicial complex $\Delta$,
the maximum likelihood estimate of $\mu$ is the point
$\hat{\mu} \in \mathcal{M}_\Delta$ such that $\hat{{\bf m}} \equiv
\exp(\hat{\mu})$ best approximates the unknown mean ${\bf m} =
\mathbb{E} ({\bf n})$ in the sense that it maximizes the probability
of observing the actual table ${\bf n}$, i.e., joint distribution of the counts
${\bf n}$ as a function of the mean  vector ${\bf m}$. This probability is also
 known as the likelihood function when we express it as a function of the
  parameters ${\bf m}$ given the data  ${\bf n}$.  The log-likelihood function 
  $\ell({\bf m})$ is the logarithm of the likelihood function.

For a given observed table ${\bf n}$, we can write  the log-likelihood as:
\[
\ell({\bf m}) = \log\mathrm{Pr} \left( {\bf n}(i) \mid{\bf m}(i), i \in d \right) 
= \sum_{i \in d} {\bf n}(i) \log{\bf m}(i) - \sum_{i \in d} {\bf m}(i) + C_{\bf n}
\]
where $C_{\bf n}$ is the logarithm of the normalization constant and
depends only on $\mathbf{n}$ and the particular sampling scheme. 
For a hierarchical model $\Delta$, we can reparametrize the log-likelihood as: 
\[
\ell(\mu) = (\mathcal{P}_{\Delta} {\bf n})^{T} \mu - \sum_{i \in d}
\exp(\mu(i)) + C_{\bf n} 
\]
where $\mathcal{P}_{\Delta}$ is the projection matrix onto $\mathcal{M}_{\Delta}$.

The {\it maximum likelihood estimate} of $\mu$ is then the vector
$\hat{\mu} \in \mathcal{M}_\Delta$ such that:
\[
\ell({\hat{\mu}}) =\sup_{\mu \in \mathcal{M}_{\Delta}} \ell(\mu)
\]
If the supremum is not attained, then the MLE is not defined.
 The log-likelihood depends on the observed table ${\bf n}$ only through
$\mathcal{P}_{\Delta} {\bf n}$ or, equivalently, since the rows of
$A_{\Delta}$ span $\mathcal{M}_{\Delta}$, the vector ${\bf t} =
A_{\Delta}{\bf n}$. Therefore, in order to establish the existence and
find the numerical value of the MLE, we need only  observe ${\bf
t}$, the vector of margins of the observed table; these are known as 
the {\it minimal sufficient statistics} for the model.

Surprisingly, the study of the conditions of existence of the MLE
has received only limited attention in the statistical literature.
Essentially all 
 available results are variations of the following theorem due to
 \citet{HAB74}: 
\begin{thm}\label{thm:hab}
 Under any of the three sampling schemes described above, a necessary and
 sufficient condition for the existence of the MLE is that there
 exists ${\bf z} \in \ker(A_\Delta)$ such that ${\bf n}+ {\bf z} >
 {\bf 0}$. 
\end{thm}

For a strengthening  of Theorem \ref{thm:hab} see \citet{GMS02}.
For a given log-linear model $\Delta$, define the marginal cone
$P_\Delta=C(A_{\Delta})$ to be the set of minimal sufficient margins, ${\bf t}$,
where, for any matrix $A$, $C(A)$ indicates the cone generated by its
columns. Let $\mathrm{relint} (P_\Delta)$ denote the relative interior
of $P_\Delta$, defined as the interior of $P_\Delta$ with respect to
its embedding into the smallest linear hull containing it. Then, the
following corollary provides a polyhedral reinterpretation of the
conditions for the existence of the MLE:
\begin{cor}\label{thm:mle}
Under any of the three sampling schemes, the MLE for the
mean vector ${\bf m}$ exists if and only if the margins ${\bf t} =
A_{\Delta} {\bf n}$ belong to $\mathrm{relint} \left( P_{\Delta}
\right)$.
\end{cor}
\begin{pf}
A vector of margins $\mathbf{t}$ lies in the relative interior of the
polyhedral 
cone $P_\Delta$ if and only if there is a table $\mathbf{x}$ with
strictly positive cells such 
that $A\mathbf{x} = \ttt$.  Theorem~\ref{thm:hab} then implies that
the MLE exists if and only if $\ttt \in {\rm relint}(P_\Delta)$. 
\qed
\end{pf}

\section{Determining the existence of the MLE}\label{sec:exist}

In this section, we describe algorithms for determining whether the
MLE for a given table $\mathbf{n}$ and model $\Delta$ exists.
To make the mathematical statements in this section
concise, we assume that 
$A_\Delta$ contains extra rows determined by the faces
of $\Delta$ in addition to those rows determined by the facets of
$\Delta$.  Since this over-parameterization does not change the row
span $\mathcal{M}_\Delta$, the matrix $A_\Delta$ describes the same
hierarchical log-linear model.  To implement the algorithms
we describe, one can relax this condition on $A_\Delta$.

By Corollary~\ref{thm:mle}, the maximum likelihood estimate
does not exist if and only if the vector of observed margins $\ttt=A_\Delta \nn$
lies on a facet of $P_\Delta$. 
Hence, we want to
show that there is a nontrivial vector $\mathbf{c}$ in the dual cone
of $P_\Delta$ which attains its maximum value at $\ttt$ but does not
attain its maximum value at some other point of $P_\Delta$.  The
existence of such a $\mathbf{c}$ implies that $\ttt$ lies on a facet of
$P_\Delta$.  However, this can be decided by determining if the
polyhedral cone 
\begin{equation}\label{eqn:lin} 
F^\Delta_\nn = \{ \mathbf{c} \mid \mathbf{c}^T A_\Delta \leq
\mathbf{1}^T \cdot \mathbf{c}^T \ttt\}
\end{equation}
contains only those vectors orthogonal to
the linear hull of $P_\Delta$.

Note that this linear system involves exponentially many inequalities
in the number of random variables $K$.  We show, however, that if
the model $\Delta$ satisfies certain nice complexity properties, the
linear system (\ref{eqn:lin}) there is an equivalent
formulation  using only
polynomially many inequalities.  Since we can solve linear programs 
in polynomial time 
\citep[e.g.,][]{SCH98}, this implies the following result:

\begin{thm} \label{thm:treewidth}
There is an algorithm for deciding the triviality of the linear
program (\ref{eqn:lin})  which runs in polynomial time in the size of
the input data and 
the number of levels of each random variable whenever the simplicial
complex $\Delta$ has bounded tree width.
\end{thm}

First, we  define all of the objects in question.

\begin{defn}
 A simplicial complex $\Delta$ is reducible if there is a
 decomposition of $\Delta$ into $(\Delta_1, S, \Delta_2)$ such that
\begin{enumerate}
\item  $\Delta_1 \cup \Delta_2  = \Delta$,
\item  $|\Delta_1| \cap |\Delta_2| = S$, and
\item  $S \in \Delta_1$ and $S \in \Delta_2$.
\end{enumerate}
Here $|\Delta_i|$ denotes the underlying set of $\Delta_i$.  A
simplicial complex is called decomposable or chordal if it is 
reducible and each of $\Delta_1$ and $\Delta_2$ are either
decomposable or a simplex.
\end{defn}

\begin{defn}
The tree width $T(\Delta)$ of a simplicial complex $\Delta$ is one
less than the size of the maximal face in the smallest decomposable
complex that contains $\Delta$.  That is,
\[
T(\Delta) = \min_{\Delta \subset \Gamma}  \max_{C \in \Gamma}  |C| - 1
\]
where the minimum runs over all decomposable  $\Gamma$ with
all faces of $\Delta$ in $\Gamma$.  A decomposable simplicial complex
$\Gamma$ that attains the minimum is called a chordal triangulation of
$\Delta$.
\end{defn}

For instance, the tree width of the $K$-cycle,
$\Delta = [12][23]\cdots[(K-1)K][1K]$, is always 2 since
a $K$-cycle does not have tree width 1 (i.e., it is not a tree), and the
simplicial complex $\Gamma = [123][134]\cdots[1(K-1)K]$ is a
decomposable complex that triangulates the $K$-cycle.  We study
the $K$-cycle in more detail in Example~\ref{ex:cycle} below.

The proof of Theorem~\ref{thm:treewidth} follows from a series of
results relating the system of linear inequalities to systems of
inequalities for chordal triangulations.  Our goal is to produce a
polyhedral cone whose triviality is equivalent to the triviality of
the cone (\ref{eqn:lin}) but whose description involves fewer linear
equations and inequalities. 

\begin{lem}\label{lem:cap}
Suppose that $\Gamma$ is a model with $\Delta \subseteq \Gamma$. Then
\[
F^\Delta_\nn = 
\pi(F^{\Gamma}_\nn  \cap \{ \mathbf{c} \mid \mathbf{c}^F = 
0 \text{ with } F \in \Gamma \setminus \Delta \}),
\]
where $\pi$ is the coordinate projection of $F^{\Gamma}_{\mathbf{n}}$
to the ambient space of $F^\Delta_\nn$.  The notation $\mathbf{c}^F$
denotes the part of the vector $\mathbf{c}$ which is naturally
labeled by the face $F \in \Gamma$.
\end{lem}
\begin{pf}
By definition. \qed
\end{pf}

Suppose that $\Delta$ is reducible, with decomposition $(\Delta_1,
S,\Delta_2)$.  From the vector $\nn$ we can compute the margins with
respect to $|\Delta_1|$ and $|\Delta_2|$, which we denote by 
$\nn_1$ and $\nn_2$.

\begin{lem}\label{lem:mink}
Suppose that $\Delta$ is reducible, with decomposition $(\Delta_1, S,
\Delta_2)$ and let $\nn$ be a table.  Then
\[
F^\Delta_\nn = \iota_1(F^{\Delta_1}_{\nn_1}) +\iota_2(F^{\Delta_2}_{\nn_2})
\]
where the ``+"  indicates the Minkowski addition of the two
cones and $\iota_1$, $\iota_2$ are the natural embeddings of
$F^{\Delta_1}_{\nn_1}$ and $F^{\Delta_2}_{\nn_2}$ into the ambient
space of $F^\Delta_\nn$.
\end{lem}
\begin{pf}
Modulo the lineality space of $F^\Delta_\nn$, the extreme rays of
$F^\Delta_\nn$ are precisely the facet defining inequalities of
$P_\Delta$ on which $\ttt$ lies.  To show the claim, it suffices to
show that every facet of $P_\Delta$ comes from a facet of $P_{\Delta_1}$
or $P_{\Delta_2}$, in the sense that ${\rm dual}(P_\Delta) =
\iota_1({\rm dual}(P_{\Delta_1})) + \iota_2({\rm
dual}(P_{\Delta_2}))$.  But this amounts to showing 
that we can decide  
the  consistency of
margins for a reducible model  by checking consistency
for both component models, $\Delta_1$ and $\Delta_2$.  Now if
the margins $\ttt_1$ and $\ttt_2$ are consistent with respect to
$\Delta_1$ and $\Delta_2$ respectively, there are tables $\nn_1$ and
$\nn_2$ such that $A_{\Delta_1} \nn_1 = \ttt_1$ and $A_{\Delta_2} \nn_2
= \ttt_2$.  Then $\nn_1$ and $\nn_2$ are margins of the decomposable
model $\Delta^* = [ |\Delta_1|][|\Delta_2|]$ which satisfy the linear
consistency relation that their $S = |\Delta_1| \cap |\Delta_2|$
margins agree.  Thus, $\mathbf{t}$ are consistent $\Delta$ marginals by
\citet{LAU96}.  This completes the proof. \qed 
\end{pf}
The description of $F^\Delta_\nn$ as a Minkowski sum in Lemma~\ref{lem:mink}
does not give a description of $F^\Delta_\nn$ that is
short in terms of having few facets.  The key to
such a short description is to recall that the Minkowski sum of two
polyhedra $P + Q$, is the image of $P \times Q$ under the map $\pi$ that
sends $(x,y)$ to $x + y$.  In particular, various properties of $P + Q$
can be determined by studying properties of $P \times Q$.  If
$P$ has $m$ facets and $Q$ has $n$ facets, then $P \times Q$ has only
$m + n$ facets.   This implies that if $P$ and $Q$ have short
descriptions in terms of few facets, then so does $P \times Q$.
Lastly, linear conditions on $P + Q$ lift to linear
conditions on $P \times Q$.  Thus we can  decide if $(P + Q) \cap L$ is empty 
be considering  $(P \times Q) \cap L'$ where $L' = \pi^{-1}(L)$.
If we accumulate all of these ideas,
together with the preceding lemmas, we get the following explicit
version of Theorem~\ref{thm:treewidth}.

\begin{thm}\label{thm:prog}
Let $\Delta$ be a simplicial complex and $\Gamma$ a chordal
triangulation of $\Delta$, with facets $\Gamma_1, \ldots, \Gamma_s$.
Denote by $\nn_t$ the $\Gamma_t$ margin of $\nn$.
Then the polyhedron  $F^\Delta_\nn$ is equal to the orthogonal
complement of the linear hull of $P_\Delta$ if
and only if the  polyhedron
\begin{equation}\label{eqn:big}(F^{\Gamma_1}_{\nn_1} \times
  \cdots  \times   F^{\Gamma_s}_{\nn_s}) \bigcap \{ (\mathbf{c}_1,
  \ldots, \mathbf{c}_s) \mid
  \sum_{i = 1}^d \mathbf{c}_t^F = \mathbf{0} \text{ for all } F \in
  \Gamma \setminus \Delta \}
\end{equation}
is a linear space.  Furthermore, if $\Delta$ has bounded
tree width, the description of \ref{eqn:big} in terms of inequalities
and equations has size that is polynomial in the number of levels of
each random variable, the number of random variables and the bit
complexity of $\nn$.  The dimension of the ambient space of the 
set in \ref{eqn:big} has size polynomial in the input.
\end{thm}
\begin{pf}
This is straightforward once we unravel all of the
definitions.  The main point is that (\ref{eqn:big}) projects, under the
``Minkowski summation'' map, onto $F^\Delta_\nn$.  This is because the
set on the left of the $\cap$ projects onto $F^{\Gamma}_\nn$ and
the set on the right of the $\cap$ is the pullback of the linear
conditions which are forced in Lemma~\ref{lem:cap}.

The statement about the complexity of the description of
(\ref{eqn:big}) follows from the fact that each of the sets $F^{\Gamma_t}_{\nn_t}$ has
a description  in terms of polynomially many
facets since the cardinality of $|\Gamma_t|$ is bounded.  The number
of inequalities 
needed to describe the object on the left hand side of the $\cap$ is
just the union of $t$ (which is a polynomial in the number of random
variables) sets of inequalities which is each only polynomial in
size.  There are only polynomially many linear conditions on the right
hand side of the $\cap$ since, if the tree width of $\Delta$ is
bounded, the cardinality of $\Gamma \setminus \Delta$ is at worst
polynomial in the number of random variables.  The dimension of the
ambient space of (\ref{eqn:big}) is polynomial in the data since the
cardinality of  $|\Gamma_t|$ is bounded.  This completes the proof of
the main theorem. 
\qed
\end{pf}

\begin{exmp}[$5$-cycle]\label{ex:cycle}
Now we will describe our construction in the special case where $K=5$ and 
$\Delta$ is the $5$-cycle.   Let $\Delta =
[12][23][34][45][15]$ and let $\Gamma = [123][134][145]$ be a chordal
triangulation.  Clearly, $\Delta$ has tree width $2$ as we previously
stated.  Now we construct the system of inequalities and
equations  in Theorem~\ref{thm:prog} for $\Delta$
with respect to $\Gamma$.

The three facets of $\Gamma$ are $\Gamma_1 = [123]$,
$\Gamma_2 = [134]$, and $\Gamma_3 = [145]$.  
>From the data, we compute the matrices
$A_{\Gamma_t}$.  we determine each of the cones
$F^{\Gamma_t}_{\nn_t}$ by the polynomially many
inequalities given by
\begin{equation} \label{eqn:red}F^{\Gamma_t}_{\nn_t} = \{\mathbf{c}_t \mid \mathbf{c}_t^T
A_{\Gamma_t} \leq \mathbf{1}^T \cdot  \mathbf{c}_t^T
A_{\Gamma_t} \nn_t \}.
\end{equation}

For each $t$, the vector $\mathbf{c}_t$ divides into
blocks, one for each face $F$ of $\Gamma_t$.  Thus, when $\Gamma_{t_1}$
and $\Gamma_{t_2}$ have a nontrivial overlap, there will be some blocks, 
$\mathbf{c}_{t_1}$ and $\mathbf{c}_{t_2}$, labeled by
the same faces.  For instance, $\Gamma_1$ and $\Gamma_2$ intersect in
the face $[13]$.

The conjunction of all the inequalities in (\ref{eqn:red}) gives all
the inequalities from the description in (\ref{eqn:big}).  To deduce the
equations, we must set to zero all of the $\mathbf{c}_t$ block
corresponding faces of $\Gamma$ that are not in $\Delta$ \emph{after}
the projection.  This amounts to adding the five sets of equations:
\[
\mathbf{c}^{[123]}_{1} = \mathbf{0}, \mathbf{c}^{[134]}_2 =
\mathbf{0}, \mathbf{c}^{[145]}_3 = \mathbf{0},
\]
\[\mathbf{c}^{[13]}_1 + \mathbf{c}^{[13]}_2 = \mathbf{0}, \mbox{ and }
\mathbf{c}^{[14]}_2 + \mathbf{c}^{[14]}_3 = \mathbf{0}.
\]

Alltold, we have a system of $O(D^3)$ inequalities
and equations, where $D = \max \{ d_1, \ldots, d_5 \}$, 
to decide if the cone is a linear space (as opposed to
$O(D^5)$ in the standard representation).

\end{exmp}

\section{Three-way tables}\label{sec:three}

\subsection{Collapsing}
In this section, we let $\Delta$ be the simplicial
complex $[12][13][23]$ on three random variables with levels $p, q,
r$, corresponding to the 
log-linear  model of no three-factor effect (also referred to as no
second-order interaction).   
This is the hierarchical log-linear model on the
fewest number of random variables where the  facet structure of the
marginal cone is not
completely understood.  From a practical standpoint, the linear
programming based algorithm from Section 3 runs in polynomial time to
determine whether or not the MLE exists for a given table under  the
no three-factor 
effect model.  However, having an understanding of the facet
structure of the marginal cone provides insight into the different
possible ways that the MLE might not exist.  Even in this small
hierarchical model, the marginal cone is quite complicated.   

Denote by $P_\Delta^{p,q,r}= P_\Delta$ the
marginal cone for this model.  We now place special emphasis on
the levels and we seek to understand the
combinatorial structure of the set of facets of $P_\Delta^{p,q,r}$.
Our main tool is collapsing the $p\times q \times r$ table to a table
with fewer levels through the combination of levels.

An {\it elementary}  collapsing of $P^\bd_\Delta$ is a linear transformation
$\pi \colon P_\Delta^{\bd} \to P_\Delta^{\bd'}$ which is obtained by
replacing some random variable $X_j$ and a set $S$ of states of $X_j$
by a new random variable $X'_j$ with $d_j - |S| + 1$ states where all
the states in $S$ are mapped to a single state.  A collapsing is any
linear map $\pi 
\colon P_\Delta^{\bd} \to P_\Delta^{\bd'}$ obtained by a sequence of
elementary collapsings.  Collapsing occurs naturally in applications
where one wishes to make coarser distinctions on the states of random
variables.  For instance, a random variable which represents the
height of individuals might be collapsed to the binary random variable
whose two states are ``tall'' and ``short''.

Since a collapsing $\pi$ maps $P^\bd_\Delta$
onto $P_\Delta^{\bd'}$, for any facet $F'$ of $P_\Delta^{\bd'}$, $F =
\pi^{-1}(F')$ is a face of $P_\Delta^{\bd}$.  If $F$ is a
facet of $P_\Delta^{\bd}$, we say that $F$ is obtained by collapsing
the $d_1 \times \dots \times d_n$ table to a $d_1' \times \dots \times
d_n'$ table.  As an example of this construction, we use collapsing to
derive exponential lower bounds on the number of facets of the
marginal cone of the no three-factor effect model.

\begin{prop}\label{prop:expmany}
The number of facets of $P_\Delta^{p,q,r}$
is at least
\[
\frac 1 2 (2^p-2)(2^q-2)(2^r-2) + pq + qr + pr.
\]
\end{prop}
\begin{pf}
Up to  symmetry, the facets of a $2\times 2 \times 2$ table
are given by the conditions:
\[
\begin{smallmatrix}
  0 & \ast \\  0 & \ast
\end{smallmatrix}
\mid
\begin{smallmatrix}
  \ast & \ast \\ \ast & \ast
\end{smallmatrix}
\;\;\;\;\;\;\; \mathrm{or} \;\;\;\;\;\;\; 
\begin{smallmatrix}
  0 & \ast \\  \ast & \ast
\end{smallmatrix}
\mid
\begin{smallmatrix}
  \ast & \ast \\ \ast & 0
\end{smallmatrix}
\]
The $0/\ast$ notation means that the facet is given by the
conditions that the ``0'' entries in the table are zero and the
$\ast$ entries are non-negative.  That is, the facet described by
a $0/\ast$ pattern is the cone over the extreme rays of the marginal
cone which are marked with a $\ast$.

The first condition says that one entry in one of the margins is zero.  There
are $pq + qr + pr$ margins for a $p\times q \times r$ table.
For the second condition, any $p\times q \times r$ table can be
collapsed to a $2 \times 2 \times 2$ table in
$(2^{p-1}-1)(2^{q-1}-1)(2^{r-1} -1)$ ways.  Each of these collapsings
gives a distinct face of $P_\Delta^{p,q,r}$ of the second type in 4
different ways.   We now show that this face is in fact a facet.
For this it suffices that the dimension of the linear span of the
extreme rays of  $P^{p,q,r}_\Delta$ that are contained in this face
has dimension 
one less than the dimension of the marginal cone.  This in turn will be
implied by showing that the linear span of these extreme rays
together with any other other extreme ray not in the face contains the
entire marginal cone  $P^{p,q,r}_\Delta$.  
Without loss of generality, by applying the natural symmetry of
this problem, it follows that the extreme rays not contained in
the face $F$ are those that have indices (i.e., positions in the $p
\times q \times r $ array) in the set
\[
I =  \{ (i_1, i_2, i_3) \mid i_1 \leq k_1, i_2 \leq k_2, i_3 \leq k_3 \}
\cup \{(i_1, i_2, i_3) \mid i_1 > k_1, i_2 > k_2, i_3 > k_3 \},
\]
for some fixed values $k_1, k_2,$ and $k_3$.
We denote the extreme ray indexed by $(i_1, i_2, i_3)$ by
$e_{i_1i_2i_3}$.  Without loss of generality, we may take $e_{111}$ to
be the extreme ray not contained in $F$, by again applying the
symmetry of the cone.  Then for any index $(j_1, j_2, j_3)$ with
$j_i > k_i$ for $i=1,2,3$,
we have the relation
\[
e_{111} + e_{1j_2j_3} + e_{j_11j_3} + e_{j_1j_21} - e_{11j_3} -
e_{1j_21}-e_{j_111} = e_{j_1j_2j_3}.
\]
Since all the extreme rays on the left hand side are
contained in $F \cup \{e_{111}\}$, this implies that $e_{j_1j_2j_3}$
is contained in the linear span of $F \cup \{e_{111}\}$.  By symmetry,
all the extreme rays indexed by elements of $I$ are contained in the
linear span of $F \cup \{e_{111}\}$.  This completes the proof that
$F$ is a facet.
\qed
\end{pf}

The cones $P_\Delta^{p,q,r}$  appear in other guises in the mathematical literature. 
for example,  \citet{vlach} studied conditions
for the non-emptiness of the three-dimensional transportation
polytopes.  A three-dimensional transportation polytope is a set of
tables
$$ P_\mathbf{t} = \{\mathbf{x} \in \mathbb{R}^d_{\geq 0} \mid
A^{p,q,r}_\Delta \mathbf{x} = \mathbf{t} \},$$
which is nonempty if and only if $\mathbf{t} \in P^{p,q,r}_\Delta$.
Hence, his results can be reinterpreted in our language.  One such
result is:

\begin{prop}\label{prop:2pq}
All facets of $P_\Delta^{2,q,r}$  are obtained by collapsing to
$P_\Delta^{2,2,2}$.
\end{prop}

Notice that Propositions~\ref{prop:expmany} and \ref{prop:2pq}
combine to show that there are exactly
$(2^q-2)(2^r-2) + 2(q+r) + qr$ facets of $P_\Delta^{2,q,r}$.

\subsection{Computations}
The polyhedron $P_\Delta^{p,q,r}$ is given by the positive hull of the
columns of $A_\Delta$ as a cone with $pqr$ extreme rays in $\Rset^{pq
+ pr + qr}$.  Some of the rows of $A_\Delta$ are redundant: the
cone is $pq + pr + qr - p - q -r+ 1$ dimensional.  It is
generally a difficult computational problem to take convex/positive hulls in a
high dimensional space.  The best algorithms for computing the convex
hull of $n$ points in $\Rset^d$ take $O(n^{\lfloor d/2\rfloor})$ time.
Using the software {\tt polymake} by \citet{polymake}
we have computed the facets for a number of
examples. 

 The group $S_p \times S_q \times S_r$ provides a natural
action on the set of facets
of $P_{p,q,r}$ given by permuting the levels of
each random variable.  After computing all the facets, we computed
orbits under this action, which gives a better picture of the set of
facets.  The results of our computations are displayed in Table 1.

It is an interesting computational problem to use this very large
symmetry group to better compute the convex hull.  The set of symmetry
classes of facets is small, and many of these classes come from
collapsing from a smaller table.  Thus many of the facets are known
``for free'' and this information should be used to compute the other
facets.  Also, the symmetry group is transitive on the extreme rays of
the cone, so in principle one could hope to compute all the facets
incident to a single extreme ray, and then use symmetry to recover the
entire cone.

Given Proposition~\ref{prop:2pq},  a natural conjecture is that
all facets are obtained by collapsing to binary tables. Unfortunately,
our computations show that the situation is remarkably more
complicated, and not all
facets of $P_\Delta^{p,q,r}$ for general $p,q,r$ are obtained by
collapsing.

\begin{exmp}[A non-collapsible facet]
The following is a facet of $P_\Delta^{4,4,4}$
that does not arise from collapsing to any smaller table.
\begin{center}
\begin{tabular}{cccc|cccc|cccc|cccc}
0&0&0&\as       &     \as&0&0&\as   &  \as&\as&0&\as   & \as&\as&\as&\as\\
0&0&\as&\as     &     \as&0&\as&\as &  \as&\as&\as&\as & 0&0&\as&0\\
0&\as&\as&\as   &    \as&\as&\as&\as&  0&\as&0&0       & 0&\as&\as&0\\
\as&\as&\as&\as &     \as&0&0&0     &  \as&\as&0&0 & \as&\as&\as&0\\
\end{tabular}
\end{center}
This example was found after examining the 39 symmetry classes of
facets of $P_\Delta^{4,4,4}$. 
\end{exmp}
\vspace{11pt}
\begin{table}
\caption{Summary of Computations.  The column ``Orbits'' counts the number
of $S_p \times S_q \times S_r$ orbits of facet types.
The column ``Collapsing'' shows the smallest table such that all facets of
$P_\Delta^{p,q,r}$ are obtained by collapsing to it. }\label{tab:comp}
\centering
\begin{tabular}{|ccc|ccccc|}
\hline
p & q & r & Dim & Extreme rays  & Facets    & Orbits & Collapsing\\
\hline
2 & 2 & 2 & 7  & 8  & 16 & 4 & 2 2 2\\
2 & 2 & 3 & 10  & 12 & 28 & 4 & 2 2 2\\
2 & 2 & 4 & 13 & 16 & 48 & 5 & 2 2 2\\
2 & 3 & 3 & 14 & 18 & 57 & 5 & 2 2 2\\
2 & 3 & 4 & 18 & 24 & 110& 6 & 2 2 2 \\
\hline
3 & 3 & 3 & 19 & 27 & 207 & 8 & 3 3 3\\
3 & 3 & 4 & 24 & 36 & 717 & 10 &3 3 3\\
3 & 3 & 5 & 29 & 45 & 2379 & 13 &3 3 3\\
3 & 3 & 6 & 34 & 54 & 7641 & 17 &3 3 3\\
3 & 3 & 7 & 39 & 63 & 23991 & 20 &3 3 3\\
\hline
3 & 4 & 4 & 30 & 48 & 4948 & 16 &3 4 4 \\
3 & 4 & 5 & 36 & 60 & 29387 & 24 & 3 4 4\\
3 & 4 & 6 & 42 & 72 & 153858 & 35 & 3 4 4\\
\hline
3 & 5 & 5 & 43 & 75 & 306955 & 42 & 3 5 5\\
\hline
4 & 4 & 4 & 37 & 64 & 113740 & 39 & 4 4 4\\
\hline
\end{tabular}
\end{table}

Based on our computations  (see Table~\ref{tab:comp}), we are led to
the following conjecture.
\begin{conj}\label{conj:col}
Suppose that $p \leq q \leq r$.  Then all facets of $P_\Delta^{p,q,r}$ are
obtained by collapsing from facets of $P_\Delta^{p,q,q}$.
\end{conj}

In general, it is true that if we fix $p$ and $q$, there exists an $r$
such that for all $r ' \geq r$, all facets of $P_\Delta^{p,q,r'}$ are
obtained by collapsing from facets of $P_\Delta^{p,q,r}$.  This
follows by noting that in a facet not obtained by collapsing, no two
slices can have the same $0/\ast$ pattern.  Since for fixed $p$ and
$q$ there are only finitely many patterns, the statement follows.
Conjecture \ref{conj:col} merely asserts that the minimal such $r$ is $q$.

\section{Summary}

We have given a polyhedral description of  the statistical problem of
determining the existence or nonexistence of the maximum likelihood
estimate for a hierarchical log-linear model for a multi-way
contingency table.  The computational implementation of this
description in principle allows statisticians to
explore for the first time the implication of patterns of zeros in
large sparse tables that lead to nonexistence and thus to recast the
estimation problem in terms of extended log-linear models for a
corresponding incomplete contingency table (c.f., \cite{HAB74}).  
There are further ties to this extended estimation problem  inherent
in the algebraic geometry description of log-linear models in terms of
Gr\"obner bases given by \cite{GMS02}.

\section*{Acknowledgments}

Nicholas Eriksson was supported by an NDSEG fellowship.
Stephen Fienberg and Alessandro Rinaldo were supported in part by National Science Foundation Grant No. EIA-0131884 to the National Institute of Statistical Sciences and Stephen Fienberg was also supported  by the Centre de Recherche en Economie et
Statistique of the  Institut National de la Statistique et des
\'{E}tudes \'{E}conomiques, Paris, France.

\bibliographystyle{elsart-harv}

\end{document}